\documentclass{amsart}
\usepackage{newlfont}
\usepackage[dvipdf]{graphicx}
\usepackage[centertags]{amsmath}
\usepackage{amsfonts}
\usepackage{amssymb}
\usepackage{amsthm}
\usepackage{caption}

\newtheorem{theorem}{theorem}[section]
\newtheorem{corollary}[theorem]{Corollary}
\newtheorem{lemma}[theorem]{Lemma}

\theoremstyle{definition}
\newtheorem{definition}[theorem]{Definition}
\theoremstyle{remark}

\newcommand{\id}{\mathrm{id}}

\title[Spheres and circles with respect to an indefinite metric]{Spheres and circles with respect to an indefinite metric on a Riemanian manifold with a skew-circulant structure}
\author{Georgi Dzhelepov, Iva Dokuzova and Dimitar Razpopov}

\begin{document}

 %%% Insert a brief summary (50-150 words)
\begin{abstract}
 We study hyper-spheres, spheres and circles, with respect to an indefinite metric, in a tangent space on a 4-dimensional differentiable manifold. The manifold is equip\-ped with a positive definite metric and an additional tensor structure of type $(1, 1)$. The fourth power of the additional structure is minus identity and its components form a skew-circulant matrix in some local coordinate system. The both structures are compatible and they determine an associated indefinite metric on the manifold.
\end{abstract}

\maketitle
\section{Introduction}

 There are various applications of the correspondences between circles and ellipses (circles and hyperbolas, circles and parabolas), as well as between spheres and other quadratic surfaces, for example in geometry, mechanics, astrophysics. Circles and spheres could be determined with respect to an indefinite metric and then their images could be obtained in Euclidean space. In this vein, we consider a circle determined with respect to an associated indefinite metric on a Riemannian manifold and the corresponding quadratic curve in Euclidean space. Also, we study a sphere (a hypersphere) determined with respect to an associated indefinite metric on a Riemannian manifold and the corresponding quadratic surface (hyper-surface) in Euclidean space.
  We will mention some papers which concern models of hyper-spheres, spheres and circles with respect to some indefinite metrics and their relations with the corresponding quadratic geometrical objects (\cite{Abe, dok-dzhe, ikawa, konderak, lopez}).

 The Hermitian manifolds form a class of manifolds with an integrable almost complex structure $J$ (\cite{GrHer}). One subclass consists of the so-called locally conformal K\"{a}hler manifolds, determined  by a special property of the covariant derivative of $J$. Some of the recent investigations of locally conformal K\"{a}hler manifolds are made in \cite{angella, cherevko, huang, moroianu-ornea, prvanovic, vilcu}.

We consider a 4-dimensional Riemannian manifold $M$, endowed with a positive definite metric $g$ and an endomorphism $S$ in a tangent space $T_{p}M$ at an arbitrary point $p$ on $M$. The fourth power of $S$ is minus identity and the components of $S$ form a skew-circulant matrix with respect to some basis of $T_{p}M$. It is supposed that $S$ is compatible with $g$. Such a manifold $(M, g, S)$ is defined in \cite{dok-raz}. In \cite{raz-dok} it is proved that $(M, g, J)$, where $J=S^{2}$, is a locally conformal K\"{a}hler  manifold.
We consider the associated metric $\tilde{g}$ on $(M, g, S)$ defined by both structures $g$ and $S$. The metric $\tilde{g}$ is necessarily indefinite and it determines space-like vectors, isotropic vectors and time-like vectors in every $T_{p}M$. We study hyper-spheres in $T_{p}M$, spheres and circles in some special subspaces of $T_{p}M$ with respect to $\tilde{g}$.

The paper is organized as follows. In Sect.~\ref{sec:2}, we recall some necessary facts, definitions and statements about the manifold $(M, g, S)$ obtained in \cite{dok-raz} and \cite{raz-dok}. In Sect.~\ref{sec:3}, we find the equation of a central hyper-sphere in $T_{p}M$ with respect to the associated metric $\tilde{g}$. In Sect.~\ref{sec:4}, we consider spheres with respect to $\tilde{g}$ in special 3-dimensional subspaces of $T_{p}M$ and obtain their equations.  In Sect.~\ref{sec:5}, we consider some special 2-planes in $T_{p}M$ and we get the equations of circles with respect to $\tilde{g}$ in these 2-planes. We interpret all equations of the curves and surfaces, studied in Sect.~\ref{sec:3}, Sect.~\ref{sec:4} and Sect.~\ref{sec:5}, in terms of $g$.

\section{Preliminaries}
\label{sec:2}
The skew-circulant matrices are Toeplitz matrices, which are well-studied in \cite{davis} and \cite{grayR}. In our work we consider a tensor structure on a 4-dimensional differentiable manifold, whose component matrix is skew-circulant. Therefore we recall the following definition.
{\em A real skew-circulant matrix\/} with the first row $(a_{1}, a_{2}, a_{3}, a_{4})\in R^{4}$ is a square matrix of the form
\[ \left(
  \begin{array}{cccc}
      a_{1} & a_{2} & a_{3} & a_{4}\\
      -a_{4} & a_{1} & a_{2} & a_{3} \\
     -a_{3} & -a_{4} & a_{1} & a_{2}\\
      -a_{2} & -a_{3} & -a_{4} & a_{1}\\
      \end{array}
\right).\]

We now introduce a manifold $(M, g, S)$ in detail. Let $M$ be a $4$-dimensional Riemannian manifold equipped with a tensor $S$ of type $(1,1)$. Let the components of $S$ form the following skew-circulant matrix in a local coordinate system:
\[ S=\left(
  \begin{array}{cccc}
      0 & 1 & 0 & 0\\
      0 & 0 & 1 & 0 \\
      0 & 0 & 0 & 1\\
      -1 & 0 & 0 & 0\\
     \end{array}
\right).\]
Then $S$ has the property
\begin{equation}\label{q4}
  S^{4}=-\id.
\end{equation}
We assume that $g$ is a positive definite metric on $M$, which satisfies the equality
\begin{equation}\label{2.12}
     g(Su, Sv)=g(u,v), \quad u, v\in \mathfrak{X}M.
\end{equation}
Such a manifold $(M, g, S)$ is introduced in \cite{dok-raz}.
The manifold $(M, g, J)$, where $J=S^{2}$, is a locally conformal K\"{a}hler  manifold (Theorem 5.3 in \cite{raz-dok}).

The associated metric $\tilde{g}$ on $(M, g, S)$, defined by
    \begin{equation}\label{defF}
   \tilde{g}(u, v)=g(u, Sv)+g(Su, v),
   \end{equation}
 is necessarily indefinite.
 Consequently, having in mind (\ref{defF}), for an arbitrary vector $v$ it is valid:
 \begin{equation}\label{s}
\tilde{g}(v, v)=2g(v, Sv)=a,\qquad a\in \mathbb{R}.
\end{equation}

According to the physical terminology we give the following
\begin{definition}\label{D1}
Let $\tilde{g}$ be the associated metric on $(M, g, S)$. If a vector $u$ satisfies $\tilde{g}(u, u)>0$ (resp. $\tilde{g}(u, u)<0$), then $u$ is a space-like (resp. a time-like) vector.
If $u$ is nonzero and satisfies $\tilde{g}(u, u)=0$, then $u$ is an isotropic vector.
\end{definition}

%Further $u$ and $v$ will stand for arbitrary vectors in the tangent space $T_{p}M$.

It is well-known that the norm of every vector $u$ of the tangent space $T_{p}M$ and the cosine of the angle between two nonzero vectors $u$ and $v$ of $T_{p}M$ are given by
   \begin{equation}\label{size}
   \|u\|=\sqrt{g(u, u)},
   \end{equation}
  \begin{equation}\label{cos}
  \cos\angle(u, v)=\frac{g(u, v)}{\|u\| \|v\|}.
  \end{equation}

A basis of type $\{u, Su, S^{2}u, S^{3}u\}$ of $T_{p}M$ is called an $S$-\textit{basis}. In this case we say that {\em the vector $u$ induces an $S$-basis of \/} $T_{p}M$. In \cite{dok-raz} the following assertions are proved.
If a vector $u$ induces an $S$-basis, then

(i) the angles between the basis vectors are
\begin{eqnarray}\label{ugli}\nonumber
& \angle(u,Su)=\angle(Su,S^{2}u)=\angle(S^{2}u, S^{3}u)=\pi-\angle(S^{3}u,u),\\
& \angle(u,S^{2}u)=\angle(Su,S^{3}u)=\frac{\pi}{2}.
\end{eqnarray}

(ii) the angle $\varphi$, determined by
\begin{equation}\label{oznachenia}
\varphi=\angle(u,Su),
\end{equation}
satisfies inequalities
\begin{equation}\label{inequalities}
  \frac{\pi}{4}<\varphi<\frac{3\pi}{4}.
\end{equation}

Next we have
\begin{theorem}\label{thmR}
Let $\tilde{g}$ be the associated metric on $(M, g, S)$ and let the vector $u$ induce an $S$-basis. The following
propositions hold true.
\begin{itemize}
  \item[(i)] Vector $u$ is space-like if and only if $\varphi\in \big(\frac{\pi}{4},\frac{\pi}{2}\big)$.
  \item[(ii)] Vector $u$ is isotropic if and only if $\varphi=\frac{\pi}{2}$.
  \item[(iii)] Vector $u$ is time-like if and only if $\varphi\in(\frac{\pi}{2},\frac{3\pi}{4})$.
\end{itemize}
\end{theorem}
\begin{proof} Using (\ref{s}), (\ref{size}), (\ref{cos}) and (\ref{oznachenia}) we get
$\tilde{g}(u, u)=2\|u\|^{2}\cos\varphi$. Having in mind Definition~\ref{D1} and inequalities (\ref{inequalities}) the proof follows. \end{proof}

Evidently, due to (\ref{q4}), (\ref{2.12}) and (\ref{defF}), we state
\begin{corollary}
If $u$ is a space-like (isotropic or time-like) vector, then $Su$, $S^{2}u$ and $S^{3}u$ are space-like (isotropic or time-like) vectors, respectively.
\end{corollary}

In the next sections we get equations of hyper-spheres, spheres and circles with respect to $\tilde{g}$ in some subspaces of $T_{p}M$ on $(M, g, S)$. Obviously, the obtained curves and surfaces do not depend on the choice of the basis. We use orthonormal bases with respect to the metric $g$ on $(M, g, S)$ to find their equations easier.  In Section~\ref{sec:3}, we use an orthonormal $S$-\textit{basis} of $T_{p}M$. The existence of such bases is proved in \cite{dok-raz}. In Section~\ref{sec:4} and Section~\ref{sec:5}, we construct orthonormal bases of  3-dimensional subspaces of $T_{p}M$ and of 2-dimensional subspaces of $T_{p}M$ with the help of an arbitrary $S$-\textit{basis}.

\section{Hyper-spheres with respect to the associated metric}\label{sec:3}

Let $\{u, Su, S^{2}u, S^{3}u\}$ be an orthonormal $S$-\textit{basis} of $T_{p}M$ with respect to the metric $g$ on $(M, g, S)$.
  If $p_{xyzt}$ is a coordinate system such that the vectors $u$, $Su$, $S^{2}u$ and $S^{3}u$ are on the axes $p_{x}$, $p_{y}$, $p_{z}$ and $p_{t}$, respectively, then $p_{xyzt}$ is orthonormal.
The radius vector $v$ of an arbitrary point $(x, y, z, t)$ of $T_{p}M$ is expressed by the equality
\begin{equation}\label{V1}
  v=xu+ySu+zS^{2}u+tS^{3}u.
\end{equation}

A hyper-sphere $s$ centered at the origin $p$, with respect to $\tilde{g}$ on $(M, g, S)$, is defined by (\ref{s}).
We apply (\ref{V1}) into (\ref{s}), and bearing in mind that $p_{xyzt}$ is an orthonormal coordinate system and also equalities (\ref{q4}) and (\ref{2.12}), we obtain the equation of $s$ as follows:
\begin{equation}\label{hyp1}
  2(xy-xt+yz+zt)=a.
\end{equation}
Now, we transform the coordinate system $p_{xyzt}$ into $p_{x'y'z't'}$ by
\begin{equation}\label{transl}
\begin{array}{ll}
   x =\frac{1}{2}(x'-y'+z'-t'),& y = \frac{\sqrt{2}}{2}(-y'+t')\\
   z = -\frac{1}{2}(x'+y'+z'+t'),& t= \frac{\sqrt{2}}{2}(-x'+z').
   \end{array}
\end{equation}
We substitute (\ref{transl}) into (\ref{hyp1}) and it takes the form
\begin{equation}\label{hyprot}
  x'^{2}+y'^{2}-z'^{2}-t'^{2}=\frac{a}{\sqrt{2}}.
\end{equation}
Evidently, in terms of $g$, we have that (\ref{hyprot}) is an equation of a 3-dimensional hyperboloid.

Therefore, we state the following
\begin{theorem}
Let $\tilde{g}$ be the associated metric on $(M, g, S)$ and let the vector $u$ induce an orthonormal $S$-basis of $T_{p}M$. If $p_{xyzt}$ is a coordinate system such that $u\in p_{x}$, $Su\in p_{y}$, $S^{2}u\in p_{z}$, $S^{3}u\in p_{t}$, then the hyper-sphere (\ref{s}) has the equation (\ref{hyprot}) with respect to the coordinate system $p_{x'y'z't'}$, obtained by the transformation (\ref{transl}) of $p_{xyzt}$.
\end{theorem}

  \begin{corollary} Let $s$ be the 3-dimensional hyperboloid (\ref{hyprot}). The following
propositions are valid.
\begin{itemize}
  \item [i)]  Every point on $s$, where $a<0$, has a time-like radius vector;
\item [ii)] Every point on $s$, where $a=0$, has an isotropic radius vector;
\item [iii)] Every point on $s$, where $a>0$, has a space-like radius vector.
\end{itemize}
\end{corollary}
\begin{proof} According to Definition~\ref{D1} and due to (\ref{s}) the statement holds. \end{proof}

\begin{corollary}
Let $s$ be the 3-dimensional hyperboloid (\ref{hyprot}). Then the intersections $\sigma_{1}$, $\sigma_{2}$, $\sigma_{3}$ and $\sigma_{4}$ between $s$ and the coordinate planes of $p_{x'y'z't'}$, respectively, are the following surfaces:
\begin{itemize}
  \item [i)] $\sigma_{1}$, $\sigma_{2}$ are hyperboloids of two sheets and $\sigma_{3}$, $\sigma_{4}$ are hyperboloids of one sheet,  in case $a>0$;
  \item [ii)] $\sigma_{1}$, $\sigma_{2}$ are hyperboloids of one sheet and $\sigma_{3}$, $\sigma_{4}$ are hyperboloids of two sheets, in case $a<0$;
  \item [iii)] $\sigma_{1}$, $\sigma_{2}$, $\sigma_{3}$ and $\sigma_{4}$ are circular cones, in case $a=0$.
\end{itemize}
\end{corollary}
\begin{proof} Using (\ref{hyprot}) and the equation of the coordinate plane $x'=0$ we get the surface $\sigma_{1}:\ \sqrt{2}(y'^{2}-z'^{2}-t'^{2})=a,\ x'=0$. Consequently, if $a>0$, then $\sigma_{1}$ is a hyperboloid of two sheet, if $a<0$, then $\sigma_{1}$ is a hyperboloid of one sheet and if $a=0$, then $\sigma_{1}$ is a circular cone. Similarly, we consider the other cases of intersections $\sigma_{2}$, $\sigma_{3}$ and $\sigma_{4}$.
\end{proof}
\begin{theorem}
Let $\tilde{g}$ be the associated metric on $(M, g, S)$ and let the vector $u$ induce an orthonormal $S$-basis of $T_{p}M$. If $p_{xyzt}$ is a coordinate system such that $u\in p_{x}$, $Su\in p_{y}$, $S^{2}u\in p_{z}$ and $S^{3}u\in p_{t}$, then $u$, $Su$, $S^{2}u$ and $S^{3}u$ are isotropic vectors and their heads lie at the surface with equations
\begin{equation}\label{circles}
   x'^{2}+y'^{2}=\frac{1}{2},\quad
   z'^{2}+t'^{2}=\frac{1}{2},
\end{equation}
where $p_{x'y'z't'}$ is the coordinate system obtained by the transformation (\ref{transl}) of $p_{xyzt}$.
\end{theorem}

\begin{proof}
 Bearing in mind (\ref{defF}) and Definition~\ref{D1} we get that $u$, $Su$, $S^{2}u$ and $S^{3}u$ are isotropic vectors with respect to $\tilde{g}$. Therefore, their heads are on the hyper-cone (\ref{hyprot}) in case $a=0$. On the other hand, these heads lie at the unit hyper-sphere with respect to $g$.
This hyper-sphere with respect to $p_{x'y'z't'}$ has the equation \begin{equation}\label{sphere}
   x'^{2} + y'^{2} + z'^{2}+t'^{2}=1.
\end{equation}
The system of (\ref{hyprot}), where $a=0$, and (\ref{sphere}) gives the intersection of a hyper-cone with a hyper-sphere. This intersection, with respect to the coordinate system $p_{x'y'z't'}$, is represented by the equivalent system (\ref{circles}).
\end{proof}

\section{Spheres in a 3-dimensional subspace of $T_{p}M$}\label{sec:4}
Let the unit vector $u$ induce an $S$-\textit{basis} of $T_{p}M$.
   Hence $u$ induces four different pyramids spanned by the following triples $\{u, Su, S^{2}u\}$, $\{Su, S^{2}u, S^{3}u\}$,  $\{u, Su, S^{3}u\}$ and $\{u, S^{2}u, S^{3}u\}$. According to (\ref{2.12}) and (\ref{ugli}), the first and the second pyramid constructed on these basis vectors are equal, as well as the third and the fourth pyramid are also equal. Thus we will investigate only the subspaces with bases defined by the first and the third pyramid.

  \subsection{A sphere in the $3$-dimensional subspace of $T_{p}M$, spanned by vectors $u$, $Su$ and $S^{2}u$}

 \begin{lemma} Let $\alpha_{1}$ be a subspace of $T_{p}M$ with a basis $\{u, Su, S^{2}u\}$. The system of vectors $\{e_{1}, e_{2}, e_{3}\}$, determined by the equalities
 \begin{eqnarray}\label{orth-base}
   e_{1}=u,\quad
   e_{2}=\frac{(-\cos\varphi)u+Su-(\cos\varphi)S^{2}u}{\sqrt{1-2\cos^{2}\varphi}}, \quad e_{3}=S^{2}u
 \end{eqnarray}
 form an orthonormal basis of $\alpha_{1}$.
 \end{lemma}
 \begin{proof}
 Using (\ref{2.12}), (\ref{ugli}) and (\ref{orth-base}) we obtain
 $g(e_{1},e_{1})=g(e_{2},e_{2})=g(e_{3},e_{3})=1$ and $g(e_{1},e_{2})=g(e_{2},e_{3})=g(e_{1},e_{3})=0.$
 \end{proof}

The coordinate system $p_{xyz}$ such that $e_{1}\in p_{x}$, $e_{2}\in p_{y}$ and $e_{3}\in p_{z}$ is orthonormal.

A sphere $s_{1}$ in $\alpha_{1}$ centered at the origin $p$, with respect to $\tilde{g}$ on $(M, g, S)$, is defined by (\ref{s}). In the next statement we get the equation of $s_{1}$ with respect to the orthonormal coordinate system $p_{xyz}$.
 \begin{theorem}
Let $\tilde{g}$ be the associated metric on $(M,g,S)$ and let $\alpha_{1}$ be a 3-dimensional subspace of $T_{p}M$ with a basis $\{u, Su, S^{2}u\}$. If $e_{1}$, $e_{2}$ and $e_{3}$ are determined by (\ref{orth-base}) and $p_{xyz}$ is a coordinate system such that $e_{1}\in p_{x}$, $e_{2}\in p_{y}$, $e_{3}\in p_{z}$,  then the equation of the sphere $s_{1}$ in $\alpha_{1}$ is given by
\begin{equation}\label{s2}
 2(\cos\varphi)(x^{2}- y^{2}+ z^{2})+2\sqrt{1-2\cos^{2}\varphi}(xy+yz)=a.
\end{equation}
\end{theorem}
\begin{proof}
The radius vector $v$ of an arbitrary point $(x, y, z)$ on $\alpha_{1}$ is expressed by
 $ v=xe_{1}+ye_{2}+ze_{3}.$
We apply the latter equality into \eqref{s} and  we find
\begin{eqnarray}\label{s-pom}\nonumber
  \tilde{g}(e_{1},e_{1})x^{2}+\tilde{g}(e_{2},e_{2})y^{2}+ \tilde{g}(e_{3},e_{3})z^{2}+2\tilde{g}(e_{1},e_{2})xy\\+2\tilde{g}(e_{1},e_{3})xz+2\tilde{g}(e_{2},e_{3})yz=a.
\end{eqnarray}
Using \eqref{2.12}, \eqref{defF}, \eqref{ugli} and \eqref{orth-base}, we obtain
\begin{eqnarray*}%\label{g-ugli2}
\tilde{g}(e_{1},e_{1})=2\cos\varphi,\quad \tilde{g}(e_{3},e_{3})=2\cos\varphi,\quad \tilde{g}(e_{2},e_{2})=-2\cos\varphi,\\ \tilde{g}(e_{1},e_{2})=\tilde{g}(e_{2},e_{3})=\sqrt{1-2\cos^{2}\varphi}, \quad  \tilde{g}(e_{1},e_{3})=0.
\end{eqnarray*}
Substituting the latter equalities into \eqref{s-pom} we get \eqref{s2}.
 \end{proof}

Now, we transform the coordinate system $p_{xyz}$ into $p_{x'y'z'}$ by
\[
  x=\frac{1}{\sqrt{2}}x'+\lambda_{1}y'+\mu_{1}z',\quad y=\lambda_{2}y'+\mu_{2}z',\quad z=-\frac{1}{\sqrt{2}}x'+\lambda_{1}y'+\mu_{1}z',
\]
where
\begin{equation}\label{lambda-mu}
\begin{array}{ll}
 \lambda_{1}=\frac{1}{2}\sqrt{1+\sqrt{2}\cos\varphi},& \lambda_{2}=\frac{\sqrt{2}}{2}\sqrt{1-\sqrt{2}\cos\varphi},\\
  \mu_{1}=\frac{1}{2}\sqrt{1-\sqrt{2}\cos\varphi},& \mu_{2}=-\frac{\sqrt{2}}{2}\sqrt{1+\sqrt{2}\cos\varphi}.
  \end{array}
\end{equation}
Therefore the equation \eqref{s2} takes the form
\begin{equation}\label{tr-s}
2\cos\varphi x'^{2}+\sqrt{2}y'^{2}-\sqrt{2}z'^{2}=a.
\end{equation}

\begin{corollary} Let $s_{1}$ be the surface, determined by \eqref{tr-s} in case $a=0$. The following statements hold true.
\begin{itemize}
  \item [i)] If $\varphi\neq \frac{\pi}{2}$, then $s_{1}$ is a cone;
    \item [ii)] If $\varphi=\frac{\pi}{2}$, then $s_{1}$ separates into two planes $z'=\pm y'$.
\end{itemize}
\end{corollary}
\begin{corollary} Let $s_{1}$ be the surface, determined by \eqref{tr-s} in case $a>0$. The following statements hold true.
\begin{itemize}
  \item [i)] If  $\varphi\in(\frac{\pi}{4}, \frac{\pi}{2})$, then $s_{1}$ is a hyperboloid of one sheets;
    \item [ii)] If $\varphi\in(\frac{\pi}{2}, \frac{3\pi}{4})$, then $s_{1}$ is a hyperboloid of two sheet;
     \item [iii)] If $\varphi=\frac{\pi}{2}$, then $s_{1}$ is a hyperbolic cylinder.
\end{itemize}
\end{corollary}

\begin{corollary}\label{cor4.5} Let $s_{1}$ be the surface, determined by \eqref{tr-s} in case $a<0$. The following statements hold true.
\begin{itemize}
  \item [i)] If $\varphi\in(\frac{\pi}{4}, \frac{\pi}{2})$, then $s_{1}$ is a hyperboloid of two sheets;
    \item [ii)] If $\varphi\in(\frac{\pi}{2}, \frac{3\pi}{4})$, then $s_{1}$ is a hyperboloid of one sheet;
    \item [iii)] If $\varphi=\frac{\pi}{2}$, then $s_{1}$ is a hyperbolic cylinder.
\end{itemize}
\end{corollary}

\subsection{A sphere in the $3$-dimensional subspace of $T_{p}M$, spanned by vectors $u$, $Su$ and $S^{3}u$}

 \begin{lemma} Let $\alpha_{2}$ be a subspace of $T_{p}M$ with a basis $\{u, Su, S^{3}u\}$. The system of vectors $\{e_{1}, e_{2}, e_{3}\}$, determined by the equalities
 \begin{equation}\label{orth-base2}
   e_{1}=Su,\quad
   e_{2}=\frac{u-(\cos\varphi)Su+(\cos\varphi)S^{3}u}{\sqrt{1-2\cos^{2}\varphi}}, \quad e_{3}=S^{3}u,
 \end{equation}
 is an orthonormal basis of $\alpha_{2}$.
 \end{lemma}
 \begin{proof}
 Using \eqref{q4}, \eqref{2.12}, \eqref{ugli} and \eqref{orth-base2} we obtain
$g(e_{1}, e_{2})=g(e_{2}, e_{3})=g(e_{1}, e_{3})=0$ and  $g(e_{1}, e_{1})=g(e_{2}, e_{2})=g(e_{3}, e_{3})=1$.
 \end{proof}

The coordinate system $p_{xyz}$ such that the vectors $e_{1}$, $e_{2}$ and $e_{3}$ are on the axes $p_{x}$, $p_{y}$ and $p_{z}$, respectively, is orthonormal.

A sphere $s_{2}$ in $\alpha_{2}$ centered at the origin $p$, with respect to $\tilde{g}$ on $(M, g, S)$, is defined by (\ref{s}). In the next statement we get the equation of $s_{2}$ with respect to the orthonormal coordinate system $p_{xyz}$.
 \begin{theorem}
Let $\tilde{g}$ be the associated metric on $(M,g,S)$ and let $\alpha_{2}$ be a 3-dimensional subspace of $T_{p}M$ with a basis $\{u, Su, S^{2}u\}$. If $e_{1}$, $e_{2}$ and $e_{3}$ are determined by (\ref{orth-base2}) and $p_{xyz}$ is a coordinate system such that $e_{1}\in p_{x}$, $e_{2}\in p_{y}$, $e_{3}\in p_{z}$, then the equation of the sphere $s_{2}$ in $\alpha_{2}$ is given by
\begin{equation}\label{surf2}
 2(\cos\varphi)(x^{2}- y^{2}+ z^{2})+2\sqrt{1-2\cos^{2}\varphi}(xy-yz)=a.
\end{equation}
\end{theorem}
\begin{proof}
The radius vector $v$ of an arbitrary point $(x, y, z)$ on $\alpha_{2}$ is expressed by
  $v=xe_{1}+ye_{2}+ze_{3}.$
Then (\ref{s}) takes the form
\begin{eqnarray}\label{s-pom2}\nonumber
  \tilde{g}(e_{1},e_{1})x^{2}+\tilde{g}(e_{2},e_{2})y^{2}+ \tilde{g}(e_{3},e_{3})z^{2}+2\tilde{g}(e_{1},e_{2})xy\\+2\tilde{g}(e_{1},e_{3})xz+2\tilde{g}(e_{2},e_{3})yz=a.
\end{eqnarray}
By (\ref{2.12}),  (\ref{defF}), (\ref{ugli}) and (\ref{orth-base2}) we obtain
\begin{eqnarray*}%\label{g-ugli2}
\tilde{g}(e_{1},e_{1})=2\cos\varphi,\quad \tilde{g}(e_{3},e_{3})=2\cos\varphi,\quad \tilde{g}(e_{2},e_{2})=-2\cos\varphi,\\ \tilde{g}(e_{1},e_{2})=\sqrt{1-2\cos^{2}\varphi}, \quad  \tilde{g}(e_{1},e_{3})=0, \quad \tilde{g}(e_{2},e_{3})=-\sqrt{1-2\cos^{2}\varphi}.
\end{eqnarray*}
Substituting the latter equalities into \eqref{s-pom2} we get \eqref{surf2}.
 \end{proof}

After transformation of the coordinate system $p_{xyz}$ into $p_{x'y'z'}$ by
\[
  x=\frac{1}{\sqrt{2}}x'+\lambda_{1}y'+\mu_{1}z',\quad y=\lambda_{2}y'+\mu_{2}z',\quad z=\frac{1}{\sqrt{2}}x'-\lambda_{1}y'-\mu_{1}z'
\]
with \eqref{lambda-mu}, the equation \eqref{surf2} takes the form
\[
2\cos\varphi x'^{2}+\sqrt{2}y'^{2}-\sqrt{2}z'^{2}=a.
\]
The above equation is the same as (\ref{tr-s}).

\section{Circles in a special 2-planes of $T_{p}M$}\label{sec:5}
Let the unit vector $u$ induce an $S$-\textit{basis} of $T_{p}M$. Now we study circles in three different subspaces $\beta_{1}$, $\beta_{2}$ and $\beta_{3}$ spanned by 2-planes $\{u, S^{2}u\}$, $\{u, Su\}$ and $\{u, S^{3}u\}$, respectively.
\subsection{Circles in the 2-plane $\beta_{1}$}
Due to (\ref{ugli}) it is true that both vectors $u$ and $S^{2}u$ form an orthonormal basis of $\beta_{1}$.
We construct a coordinate system $p_{xy}$ on $\beta_{1}$, such that $u$ is on the axis $p_{x}$ and $S^{2}u$ is on the axis $p_{y}$.
Therefore $p_{xy}$ is an orthonormal coordinate system of $\beta_{1}$.

\begin{lemma} The system $\{u, S^{2}u\}$ satisfies the following equalities:
\begin{equation}\label{circle1}
  \tilde{g}(u,u)=\tilde{g}(S^{2}u,S^{2}u)=2\cos\varphi ,\quad \tilde{g}(u,S^{2}u)=0.
\end{equation}
\end{lemma}
\begin{proof}
  From (\ref{2.12}), (\ref{defF}), (\ref{ugli}) we get (\ref{circle1}) by direct calculations.
\end{proof}
A circle $k_{1}$ in $\beta_{1}$ centered at the origin $p$, with respect to $\tilde{g}$ on $(M,g,S)$, is defined by (\ref{s}). Now we obtain the equation of $k_{1}$ with respect to $p_{xy}$.
\begin{theorem}
Let $\tilde{g}$ be the associated metric on $(M, g, S)$ and let $\beta_{1}$ be a $2$-plane in $T_{p}M$ with a basis $\{u, S^{2}u\}$. If $p_{xy}$ is a coordinate system such that $u\in p_{x}$, $Su\in p_{y}$, then the equation of the circle $k_{1}$ in $\beta_{1}$ is given by
\begin{equation}\label{k-P}
 2\cos\varphi x^{2} +2\cos\varphi y^{2}=a,\quad \varphi\neq \frac{\pi}{2}.
\end{equation}
\end{theorem}

\begin{proof}
The radius vector $v$ of an arbitrary point on $\beta_{1}$ is expressed by
\begin{equation}\label{v-za-p}
  v=xu+yS^{2}u,
\end{equation}
which implies $S^{2}v=-yu+xS^{2}u$. Then from (\ref{s}), (\ref{circle1}) and (\ref{v-za-p}) it follows (\ref{k-P}).
\end{proof}
\begin{corollary} Let $k_{1}$ be the curve determined by (\ref{k-P}). Then $k_{1}$ is a circle in case when $a> 0$ and $\varphi\in \big(\frac{\pi}{4},\frac{\pi}{2}\big)$, or in case when $a< 0$ and $\varphi\in \big(\frac{\pi}{2},\frac{3\pi}{4}\big)$. The curve $k_{1}$ degenerates into the point $p$ in case $a=0$.
\end{corollary}
We note that a 2-plane $\delta=\{u, Ju\}$, where $\delta=J\delta$, is known as $J$-invariant section of $T_{p}M$ on an almost Hermitian manifold $(M, g, J)$.
Therefore, the 2-plane $\beta_{1}=\{u, S^{2}u\}$ is a $J$-invariant section of $T_{p}M$ on the manifold $(M, g, J)$, $J=S^{2}$.

\subsection{Circles in the 2-plane $\beta_{2}$}
\begin{lemma} Let $\beta_{2}$ be the $2$-plane spanned by unit vectors $u$ and $Su$. The system of vectors $\{e_{1}, e_{2}\}$, determined by the equalities
\begin{equation}\label{i-j}
 e_{1}=\frac{1}{\sqrt{2(1+\cos\varphi)}}(u+Su),\quad e_{2}=\frac{1}{\sqrt{2(1-\cos\varphi)}}(-u+Su),
\end{equation}
is an orthonormal basis of $\beta_{2}$.
\end{lemma}
\begin{proof}
Using  (\ref{ugli}) and (\ref{i-j}), we calculate $g(e_{1}, e_{2})=0$, $g(e_{1}, e_{1})=g(e_{2}, e_{2})=1$.
\end{proof}
We construct a coordinate system $p_{xy}$ on $\beta_{2}$, such that $e_{1}$ is on the axis $p_{x}$ and $e_{2}$ is on the axis $p_{y}$, i.e.  $p_{xy}$ is  orthonormal.
\begin{lemma} The system $\{e_{1}, e_{2}\}$ satisfies the following equalities:
\begin{equation}\label{circle2}
  \tilde{g}(e_{1},e_{1})=\frac{2\cos\varphi+1}{1+\cos\varphi},\quad\tilde{g}(e_{2},e_{2})=\frac{2\cos\varphi-1}{1-\cos\varphi} ,\quad \tilde{g}(e_{1},e_{2})=0.
\end{equation}
\end{lemma}
\begin{proof}
  Using (\ref{2.12}),  (\ref{defF}) and (\ref{ugli}) we get (\ref{circle2}) by direct calculations.
\end{proof}
A circle $k_{2}$ in $\beta_{2}$ centered at the origin $p$, with respect to $\tilde{g}$ on $(M, g, S)$, is defined by (\ref{s}).
Further we obtain the equation of $k_{2}$ with respect to $p_{xy}$.
\begin{theorem}
Let $\tilde{g}$ be the associated metric on $(M, g, S)$ and let $\beta_{2}=\{u, Su\}$ be a $2$-plane in $T_{p}M$ with an orthonormal basis (\ref{i-j}). If $p_{xy}$ is a coordinate system such that $e_{1}\in p_{x}$, $e_{2}\in p_{y}$, then the equation of the circle $k_{2}$ in $\beta_{2}$ is given by
\begin{equation}\label{k-q}
 \frac{2\cos\varphi+1}{1+\cos\varphi} x^{2}+\frac{2\cos\varphi-1}{1-\cos\varphi}y^{2}=a.
\end{equation}
\end{theorem}
\begin{proof}
The radius vector $v$ of an arbitrary point on $\beta_{2}$ is
 $ v=xe_{1}+ye_{2}.$
Using the latter equality, from (\ref{s}) we get
\[  \tilde{g}(v,v)=\tilde{g}(e_{1},e_{1})x^{2}+2\tilde{g}(e_{1},e_{2})xy+\tilde{g}(e_{2},e_{2})y^{2}=a.\]

Applying (\ref{circle2}) into the above equation we find (\ref{k-q}).
\end{proof}

 According to the parameters $a$ and $\varphi$ the equation (\ref{k-q}) describes different quadratic curves. All possible values of these parameters and the corresponding types of the curve (\ref{k-q}) are studied in the Table~\ref{tab:1}.

 \subsection{Circles in the 2-plane $\beta_{3}$}
\begin{lemma} Let $\beta_{3}$ be the $2$-plane spanned by unit vectors $u$ and $S^{3}u$. The system of vectors $\{e_{1}, e_{2}\}$, determined by the equalities
\begin{equation}\label{e1-e2}
 e_{1}=\frac{1}{\sqrt{2(1-\cos\varphi)}}(u+S^{3}u),\quad e_{2}=\frac{1}{\sqrt{2(1+\cos\varphi)}}(-u+S^{3}u),
\end{equation}
is an orthonormal basis of $\beta_{3}$.
\end{lemma}
\begin{proof}
Using (\ref{ugli}) and (\ref{e1-e2}), we calculate $g(e_{1}, e_{2})=0$, $g(e_{1}, e_{1})=g(e_{2}, e_{2})=1$.
\end{proof}
We construct a coordinate system $p_{xy}$ on $\beta_{3}$, such that $e_{1}\in p_{x}$ and $e_{2}\in p_{y}$, i.e.  $p_{xy}$ is orthonormal.
\begin{lemma} The system $\{e_{1}, e_{2}\}$ satisfies the following equalities:
\begin{equation}\label{circle3}
  \tilde{g}(e_{1},e_{1})=\frac{2\cos\varphi-1}{1-\cos\varphi},\quad\tilde{g}(e_{2},e_{2})=\frac{2\cos\varphi+1}{1+\cos\varphi} ,\quad \tilde{g}(e_{1},e_{2})=0.
\end{equation}
\end{lemma}
\begin{proof}
Using (\ref{q4}), (\ref{2.12}), (\ref{defF}) and (\ref{ugli}) we get (\ref{circle3}) by direct calculations.
\end{proof}
A circle $k_{3}$ in $\beta_{3}$ centered at the origin $p$, with respect to $\tilde{g}$ on $(M, g, S)$, is defined by (\ref{s}).
In the next statement we obtain the equation of $k_{3}$ with respect to $p_{xy}$.
\begin{theorem}
Let $\tilde{g}$ be the associated metric on $(M, g, S)$ and let $\beta_{3}=\{u, S^{3}u\}$ be a $2$-plane in $T_{p}M$ with an orthonormal basis (\ref{e1-e2}). If $p_{xy}$ is a coordinate system such that $e_{1}\in p_{x}$, $e_{2}\in p_{y}$, then the equation of a circle $k_{3}$ in $\beta_{3}$ is given by
\begin{equation}\label{k-q2}
 \frac{2\cos\varphi-1}{1-\cos\varphi} x^{2}+\frac{2\cos\varphi+1}{1+\cos\varphi}y^{2}=a.
\end{equation}
\end{theorem}

\begin{proof}
The radius vector $v$ of an arbitrary point on $\beta_{3}$ is
$v=xe_{1}+ye_{2}.$
Then (\ref{s}) imply
\begin{equation}\label{k-pom3}
  \tilde{g}(v,v)=\tilde{g}(e_{1},e_{1})x^{2}+2\tilde{g}(e_{1},e_{2})xy+\tilde{g}(e_{2},e_{2})y^{2}=a.
\end{equation}

We apply (\ref{circle3}) into (\ref{k-pom3}) and we find (\ref{k-q2}).
\end{proof}
The equation (\ref{k-q2}) determines curves which are the same as the obtained ones by (\ref{k-q}). They are described in the Table~\ref{tab:1}.
\begin{table}
\caption{Curves $k_{2}$ and $k_{3}$}
\label{tab:1}       % Give a unique label
\begin{tabular}{llll}
\hline\noalign{\smallskip}
 $\varphi$ & $a$ & & $k_{2}$, $k_{3}$  \\
\noalign{\smallskip}\hline
$(\frac{\pi}{4},\frac{\pi}{3})$ & $a>0$ & & an ellipse\\
 & $a=0$ & & the point $p$  \\
 & $a<0$ & & the empty set \\
 \noalign{\smallskip}\hline
 $\frac{\pi}{3}$ & $a>0$ & & two lines $x=\pm\frac{\sqrt{3a}}{2}$ \\
 &  $a=0$ & & the line x=0 \\
 &  $a<0$ & & the empty set \\
\noalign{\smallskip}\hline\noalign{\smallskip}
 $(\frac{\pi}{3}, \frac{2\pi}{3})$ & $a>0$ & & a hyperbola \\
 &  $a=0$ & & two lines $y=\pm cx$, $c$ is a constant \\
 &  $a<0$ & & a hyperbola \\
 \noalign{\smallskip}\hline
$\frac{2\pi}{3}$ &  $a>0$ & &the empty set \\
  & $a=0$ & & the line y=0 \\
  & $a<0$ & &two lines $y=\pm\frac{\sqrt{-3a}}{2}$\\
\noalign{\smallskip}\hline
$(\frac{2\pi}{3},\frac{3\pi}{4})$ & $a>0$ & & the empty set\\
 & $a=0$ & & the point $p$ \\
 & $a<0$ & & an ellipse\\
\noalign{\smallskip}\hline
\end{tabular}
\end{table}

 \bigskip

{\small\rm\baselineskip=10pt
 \baselineskip=10pt
 \qquad Georgi Dzhelepov\par
 \qquad Faculty of Economics\par
 \qquad Department of Mathematics and Informatics\par
 \qquad Agricultural University of Plovdiv\par
 \qquad 4000 Plovdiv, Bulgaria\par
 \qquad {\tt dzhelepov@abv.bg}

 \bigskip \smallskip

 \qquad Iva Dokuzova\par
 \qquad Faculty of Mathematics and Informatics\par
 \qquad Department of Algebra and Geometry\par
 \qquad University of Plovdiv Paisii Hilendarski\par
 \qquad 24 Tzar Asen, 4000 Plovdiv, Bulgaria\par
 \qquad {\tt dokuzova@uni-plovdiv.bg}

 \bigskip \smallskip

\qquad Dimitar Razpopov\par
 \qquad Faculty of Economics\par
 \qquad Department of Mathematics and Informatics\par
 \qquad Agricultural University of Plovdiv\par
 \qquad 4000 Plovdiv, Bulgaria\par
 \qquad {\tt razpopov@au-plovdiv.bg}
 }
 \end{document}